\documentstyle[amsmath]{article}

\newtheorem{th}{Theorem}[section]

\newtheorem{cor}[th]{Corollary}
\newtheorem{defn}[th]{Definition}
\newenvironment{defn-new}{\begin{defn} \em}{\end{defn}}
\newtheorem{rem}[th]{Remark}
\newenvironment{rem-new}{\begin{rem} \em}{\end{rem}}
\newtheorem{ex}[th]{Example}
\newenvironment{ex-new}{\begin{ex} \em}{\end{ex}}

\newenvironment{notation-new}{\begin{rem} \em}{\end{rem}}

\newenvironment{agr-new}{\begin{rem} \em}{\end{rem}}

\makeatletter \@addtoreset{equation}{section} \makeatother

\makeatletter \@addtoreset{figure}{section} \makeatother

\begin{document}

\begin{center}
{\Large {\bf Second Order parallel tensor on generalized f.pk-space form and
hypersurfaces of generalized f.pk-space form}}\\[0pt]

\vspace{1cm} {\bf Punam Gupta and Sanjay Kumar Singh}
\end{center}

\noindent {\bf Abstract.}
 The purpose of the present paper to study
a second order symmetric parallel tensor in generalized f.pk-space
form. Second order symmetric parallel tensor in f.pk-space form is
combination of the associated metric tensor and $1$-forms of
structure vector fields. We prove that there does not exist second
order skew-symmetric parallel tensor in f.pk-space form. We also
deduce that there is no parallel hypersurface in a generalized
f.pk-space form but there is semi-parallel hypersurfaces in a
generalized f.pk-space form.\\

\medskip


\medskip


\section{Introduction}

In 1923, Eisenhart \cite{Eisenhart-23} proved that if a positive definite
Riemannian manifold admits a second order parallel symmetric tensor other
than a constant multiple of the metric tensor, then it is reducible. In
1926, Levy \cite{Levy-26} proved that a second order parallel symmetric
non-singular (with non-vanishing determinant) tensor in a real space form is
proportional to the metric tensor. As an improvement of the result of Levy,
Sharma \cite{Sharma} proved that any second-order parallel tensor not
necessarily symmetric in a real-space form of dimension greater than $2$ is
proportional to the Riemannian metric. Later on, many authors \cite%
{Bejan,Das,Gupta-1,Mondal,Ramesh,Sharma-1,Gupta} studied
second-order parallel tensor on various spaces and obtained the
important results. In \cite{Gherib-1}, Gherib and Belkhelfa studied
the second order symmetric parallel tensors on generalized Sasakian
space form and proved that it is proportional to metric tensor.
Recently, Belkhelfa and Mahi \cite{Belkhelfa} proved that second
order symmetric parallel tensors on ${\cal S}$-space form is linear
combination of the metric tensor. and $1$-forms of structure vector
fields with constant coefficients. They also prove the existence of
semi-parallel hypersurfaces in ${\cal S}$-space form.

On the other side, Falcitelli and Pastore \cite{Fal-Pas} extend the notion
of generalized Sasakian-space form to the generalized f.pk-space form and
discussed the constancy of the $\varphi $-sectional curvature. They also
interrelate generalized f.pk-space form with generalized Sasakian and
generalized complex space-forms. Also the ${\cal S}$-space form, Sasakian
space form are the particular cases of the generalized f.pk-space form.

The purpose of this paper is to study second order parallel tensors
on generalized f.pk-space form. We prove that second order symmetric
parallel tensors on generalized f.pk-space form is linear
combination of the associtaed positive-definite metric tensor and
$1$-forms of structure vector fields. We prove the non-existence of
non-zero second order skew-symmetric parallel tensor. We also find
out that there do not exist a parallel
hypersurface of a generalized f.pk-space forms $M^{2n+s}(F_{1},F_{2},{\cal F}%
)$, tangent to the structure vector fields. Even there exist semi-parallel
hypersurface of a generalized f.pk-space forms $M^{2n+s}(F_{1},F_{2},{\cal F}%
)$, tangent to the structure vector fields.

\section{Preliminaries}

An f.pk-manifold \cite{Blair,Duggal,Gold} is a manifold $M^{2n+s}$ on which
is defined an $f$-structure \cite{Yano} (also known as $f$-structure with
complemented frames or globally framed $f$-stuctures or $f$-structures with
parallelizable kernel i.e. f.pk-structures), that is a $(1,1)$-tensor field $%
\varphi $ satisfying
\[
\varphi ^{3}+\varphi =0,
\]%
of rank $2n$ such that the subbundle $\ker \varphi $ is parallelizable. Then
there exists a global frame $\{\xi _{i}\},i=1,\ldots ,s$ for the subbundle $%
\ker \varphi $, with dual $1$-forms $\eta ^{i}$, satisfying
\begin{equation}
\varphi ^{2}=-I+\eta ^{i}\otimes \xi _{i},  \label{eq-acm-1}
\end{equation}%
\begin{equation}
\eta ^{i}(\xi _{j})=\delta _{j}^{i},\quad \varphi \,\xi _{i}=0,\quad \eta
^{i}\circ \varphi =0.  \label{eq-acm-2}
\end{equation}%
An f.pk-structure on a manifold $M^{2n+s}$ is said to be normal if the
Nijenhuis tensor field $N=[\varphi ,\varphi ]+2d\eta ^{i}\otimes \xi _{i}$
vanishes, where $[\varphi ,\varphi ]$ is the Nijenhuis torsion of $\varphi $%
. Consider a Riemannian metric $g$ on $M^{2n+s}$ associated with an
f.pk-structure $(\varphi ,\xi _{i},\eta ^{i})$ such that
\begin{equation}
g\left( \varphi X,\varphi Y\right) =g\left( X,Y\right)
-\sum\limits_{i=1}^{s}\eta ^{i}(X)\eta ^{i}\left( Y\right) ,
\label{eq-acm-5}
\end{equation}%
for any vector fields $X,Y$ on $\Gamma (TM)$. Then an f.pk-structure is
known as a metric f.pk-structure. A manifold with a metric f.pk-structure is
known as a metric f.pk-manifold.

Let $\Phi $ be the fundamental $2$-form on $M^{2n+s}$ defined by
\[
\Phi (X,Y)=g(X,\varphi Y)=-g(\varphi X,Y),
\]%
for any vector fields $X,Y$ on $\Gamma (TM)$. A metric f.pk-structure is
said to be a ${\cal K}$-structure \cite{Blair} if it is normal and the
fundamental $2$-form $\Phi $ is closed, a manifold with a ${\cal K}$%
-structure is known as ${\cal K}$-manifold. A ${\cal K}$-structure is said
to be a ${\cal S}$-structure if $d\eta ^{i}=\Phi $ for all $i\in \{1,\ldots
,s\}$, a manifold with an ${\cal S}$-structure is known as ${\cal S}$%
-manifold. A ${\cal K}$-structure is said to be a ${\cal C}$-structure if $%
d\eta ^{i}=0$ for all $i\in \{1,\ldots ,s\}$, a manifold with a ${\cal C}$%
-structure is known as ${\cal C}$-manifold. manifold. Obviously, if $s=1$, a
${\cal K}$-manifold is a quasi Sasakian manifold, a ${\cal C}$-manifold is a
cosymplectic manifold and a ${\cal S}$-manifold is a Sasakian manifold.

The Levi-Civita connection $\nabla $ \cite{Blair,Duggal} of a metric
f.pk-manifold satisfies%
\begin{eqnarray*}
2g((\nabla _{X}\varphi )Y,Z &=&3d\Phi (X,\varphi Y,\varphi Z)-3d\Phi (X,Y,Z)
\\
&&+g(N(Y,Z),\varphi X)+2d\eta ^{j}(\varphi Y,Z)\eta ^{j}(X) \\
&&-2d\eta ^{j}(\varphi Z,Y)\eta ^{j}(X)+2d\eta ^{j}(\varphi Y,X)\eta ^{j}(Z)
\\
&&-2d\eta ^{j}(\varphi Z,X)\eta ^{j}(Y).
\end{eqnarray*}%
In particular, for ${\cal S}$-manifolds \cite{Blair}, we have $\nabla
_{X}\xi _{i}=-\varphi X,\quad i=1,\ldots ,s$.

A plane section $\Pi $ in $T_{p}M$ spanned by $X$ and $\varphi X$, where $X$
is a tangent vector orthogonal to structure vector fields, is known as $%
\varphi $-section. The sectional curvature of $\varphi $-section is called a
$\varphi $-sectional curvature. A metric f.pk-manifold $(M^{2n+s},\varphi
,\xi _{i},\eta ^{i},g)$ has pointwise constant $\varphi $-sectional
curvature if at any point $p\in M^{2n+s},\quad R_{p}(X,\varphi X,X,\varphi
X) $ does not depend on the $\varphi $-section spanned by $\{X,\varphi X\}$.

\section{Generalized f.pk- space form}

\begin{defn-new}
{\rm \cite{Fal-Pas}} A generalized f.pk- space form $M^{2n+s}(F_{1},F_{2},%
{\cal F}),$ is a metric f.pk-manifold $(M^{2n+s},\varphi ,\xi _{i},\eta
^{i},g)$ which admits smooth functions $F_{1},F_{2},{\cal F}$ such that its
curvature tensor field satisfies
\begin{eqnarray}
R(X,Y)Z &=&F_{1}\left( g(\varphi X,\varphi Z)\varphi ^{2}Y-g(\varphi
Y,\varphi Z)\varphi ^{2}X\right)  \nonumber \\
&&+F_{2}\left( g(Z,\varphi Y)\varphi X-g(Z,\varphi X)\varphi Y+2g(X,\varphi
Y)\varphi Z\right)  \nonumber \\
&&+\sum\limits_{i,j=1}^{s}F_{ij}\left( \eta ^{i}(X)\eta ^{j}(Z)\varphi
^{2}Y-\eta ^{i}(Y)\eta ^{j}(Z)\varphi ^{2}X\right.  \nonumber \\
&&\left. +g(\varphi Y,\varphi Z)\eta ^{i}(X)\xi _{j}-g(\varphi X,\varphi
Z)\eta ^{i}(Y)\xi _{j}\right) ,  \label{eq-fpk}
\end{eqnarray}%
The $\varphi $-sectional curvature of a generalized f.pk-space form $%
M^{2n+s}(F_{1},F_{2},{\cal F})$ is pointwise constant, which is $%
c=F_{1}+3F_{2}$.
\end{defn-new}

\begin{rem-new}
For $s\geq 2$, an ${\cal S}$-maniolfd with pointwise constant $\varphi $%
-sectional curvature $c$ is known as ${\cal S}$-space form \cite%
{Blair,Kobayashi}. If the given structure is ${\cal S}$-structure then we
obtain a ${\cal S}$-space form by (\ref{eq-fpk}) with $F_{1}=\dfrac{c+3s}{4}%
,F_{2}=\dfrac{c-s}{4}$ and $F_{ij}=1$ for all $i,j\in \{1,\ldots ,s\}$.
\end{rem-new}

\begin{rem-new}
For $s=1$, we obtain a generalized Sasakian-space form $%
M^{2n+1}(f_{1},f_{2},f_{3})$ with $f_{1}=F_{1}$, $f_{2}=F_{2}$ and $%
f_{3}=F_{1}-F_{11}$. If the given structure is Sasakian then (\ref{eq-fpk})
holds with $F_{11}=1,F_{1}=(c+3)/4,F_{2}=(c-1)/4$ and $%
f_{3}=F_{1}-F_{11}=(c-1)/4=f_{2}$. If the given structure is Kenmotsu then (%
\ref{eq-fpk}) holds with $F_{11}=-1,F_{1}=(c-3)/4,F_{2}=(c+1)/4$ and $%
f_{3}=F_{1}-F_{11}=(c+1)/4=f_{2}$. If the given structure is cosymplectic
then (\ref{eq-fpk}) holds with $F_{11}=0,F_{1}=c/4,F_{2}=c/4$ and $%
f_{3}=F_{1}-F_{11}=c/4$.
\end{rem-new}

Let $M^{2n+s}$ be a metric f.pk-manifold with its metric tensor $g$ and
Levi-Civita connection $\nabla $. Let $R$ denote the Riemann curvature
tensor of $M^{2n+s}$. If $H$ is a $(0,2)$-tensor which is parallel with
respect to $\nabla $ then we can show easily that
\begin{equation}
H(R(X,Y)Z,W)+H(Z,R(X,Y)W)=0.  \label{eq-H}
\end{equation}

\begin{th}
\label{th-sec} A second order parallel symmetric tensor in a generalized
f.pk- space form $M^{2n+s}(F_{1},F_{2},{\cal F})$ is a linear combination of
the associated positive definite metric tensor and $1$-forms of structure
vector fields if $s\geq 2$ with assumption that ${\cal F}$ has atleast one
non-zero function.
\end{th}

\noindent{\bf Proof. }Using (\ref{eq-fpk}) in (\ref{eq-H}) and taking $Y=\xi _{k}$ $%
(k\in \{1,\ldots s\})$ and $Z=W$, we get

\begin{equation}
H\left( \sum\limits_{j=1}^{s}F_{kj}\left( g(\varphi X,\varphi Z)\xi
_{j}+\eta ^{j}(Z)\varphi ^{2}X\right) ,Z\right) =0  \label{q-1}
\end{equation}%
On solving and taking $Z=\xi _{r}(r\in \{1,\ldots s\})$, we get
\begin{equation}
H(X,\xi _{r})=\sum\limits_{\alpha ,\beta =1}^{s}\eta ^{\alpha }(X)H(\xi
_{\beta },\xi _{r}),\quad F_{kr}\not=0.  \label{eq-2}
\end{equation}%
Using (\ref{eq-2}) in (\ref{q-1}), we have
\begin{equation}
H(X,Z)=g(X,Z)\sum\limits_{\gamma =1}^{s}H(\xi _{j},\xi _{\gamma
})-\sum\limits_{\alpha ,\beta ,\gamma =1}^{s}\eta ^{\alpha }(X)\eta ^{\beta
}(Z)H(\xi _{j},\xi _{\gamma })+\sum\limits_{\alpha ,\beta ,\gamma
=1}^{s}\eta ^{\alpha }(X)\eta ^{j}(Z)H(\xi _{\beta },\xi _{\gamma }).
\label{eq-H1}
\end{equation}

\begin{rem-new}
For $s=1$, the (\ref{eq-H1}) reduces to $H(X,Z)=g(X,Z)H(\xi ,\xi )$.
\end{rem-new}

\begin{rem-new}
For an ${\cal S}$-maniolfd, $H(\xi _{j},\xi _{\gamma })$ are constant for
all $j,\gamma \in \{1,\ldots s\}$. Therefore we can say that a second order
parallel symmetric tensor in a ${\cal S}$-space form is a linear combination
of the associated positive definite metric tensor and $1$-forms of structure
vector fields with constant coefficients if $s\geq 2$ \cite{Belkhelfa}.
\end{rem-new}

\begin{rem-new}
For $s=1$ and if the manifold is Sasakian then a second order parallel
symmetric tensor in a Sasakian-space form or generalized Sasakian space form
is proportional to the associated positive definite metric tensor.
\end{rem-new}

As an application of Theorem \ref{th-sec}, now we consider the Ricci tensor
of the manifold. Since we know that the Ricci tensor $S$ of the manifold is
symmetric $(0,2)$-tensor.

\begin{defn-new}
A non-flat semi-Riemannian manifold $M$ is said to be  pseudo-Ricci
symmetric \cite{Chaki} if the Ricci tensor $S$ of type $(0,2)$ of the
manifold is non-zero and satisfies the condition
\[
(\nabla _{X}S)(Y,Z)=2\alpha (X)S(Y,Z)+\alpha (Y)S(X,Z)+\alpha (Z)S(X,Y),
\]%
for all vector fields $X,Y,Z$, where $\alpha $ is a non-zero $1$-form. If $%
\alpha =0$, then the manifold reduces to Ricci symmetric manifold.
\end{defn-new}

\begin{defn-new}
A semi-Riemannian manifold $M$ is said to be Ricci-semisymmetric \cite%
{Deszcz} if its Ricci tensor $S$ satisfies $R$\textperiodcentered $S=0$,
that is,
\[
S(R(X,Y)Z,W)+S(Z,R(X,Y)W)=0
\]%
for all vector fields $X,Y,Z,W$.
\end{defn-new}

\begin{defn-new}
A semi-Riemannian manifold $M$ is said to be Ricci flat if its Ricci tensor $%
S$ satisfies $S=0$.
\end{defn-new}

\begin{cor}
A Ricci symmetric generalized f.pk- space form $M^{2n+s}(F_{1},F_{2},{\cal F}%
)$ is a linear combination of the associated positive definite metric tensor
and $1$-forms of structure vector fields if $s\geq 2$ with assumption that $%
{\cal F}$ has atleast one non-zero function.
\end{cor}

In \cite[Theorem 2.2]{Arslan}, it was shown that if $M$ is pseudo
Ricci symmetric manifold such that ${\rm div}R=0$, then $M$ is Ricci
flat. We know that Ricci flat condition implies Ricci symmetric but
converse need not be. Therefore, we can state that

\begin{cor}
Let $M^{2n+s}(F_{1},F_{2},{\cal F})$ be a pseudo-Ricci symmetric generalized
f.pk- space form such that ${\cal F}$ has atleast one non-zero function and $%
s\geq 2$. Then the linear combination of the associated positive definite
metric tensor and $1$-forms of structure vector fields is zero.
\end{cor}

Now, we consider skew-symmetric parallel tensor on generalized f.pk- space
form.

\begin{th}
There does not exist a non-zero second order parallel skew-symmetric tensor
in a generalized f.pk- space form $M^{2n+s}(F_{1},F_{2},{\cal F})$ with
assumption that ${\cal F}$ has atleast one non-zero function.
\end{th}

\noindent {\bf Proof. }Consider (\ref{eq-H}) with $H$ as a non-zero second
order parallel skew-symmetric tensor, using (\ref{eq-fpk}) and taking $%
Y=W=\xi _{k}$ $(k\in \{1,\ldots s\})$ with assumption that ${\cal F}$ has
atleast one non-zero function, we obtain
\begin{equation}
H(Z,X)=\eta ^{k}(Z)H(X,\xi _{k})+\eta ^{k}(X)H(Z,\xi _{k}).  \label{eq-skew}
\end{equation}
Let $A$ be the dual $(1,1)$-type tensor which is metrically equivalent to $H$%
, that is, $H(X,Y)=g(AX,Y)$. Then (\ref{eq-skew}) is equivalent to
\begin{equation}
AX=-g(AX,\xi _{k})\xi _{k}+\eta ^{k}(X)A\xi _{k}.  \label{eq-skew-1}
\end{equation}%
Now, taking the inner product of (\ref{eq-skew-1}) with $\xi _{k}$, we have $%
H(X,\xi _{k})=0$. Therefore (\ref{eq-skew-1}) reduces to
\begin{equation}
AX=\eta ^{k}(X)A\xi _{k}.  \label{eq-skew-2}
\end{equation}%
Taking the inner product of (\ref{eq-skew-2}) with $Y$, we have $H(X,Y)=0$.

\section{Parallel submanifolds of f.pk-space forms}

Let $G$ be an immersed hypersurface of $M^{2n+s}(F_{1},F_{2},{\cal F})$.
Then the formulas of Gauss and Weingarten are

\[
\tilde{\nabla}_{X}Y=\nabla _{X}Y+h(X.Y)N,
\]%
\[
\tilde{\nabla}_{X}N=-SX,
\]%
where $X$ and $Y$ are tangent vector fields, $N$\ a unit normal vector
normal to $G$, $h$ the second fundamental form and $S$ the shape operator of
$G$. Note that $h$ and $S$ are related by $h(X,Y)=g(SX,Y)$. In a
hypersurface, the $(0,4)$ tensor field $\tilde{R}.h$ is defined by
\begin{equation}
\tilde{R}.h(X,Y,Z,W)=-h(\tilde{R}(X,Y)Z,W)-h(Z,\tilde{R}(X,Y)W)
\label{eq-hy}
\end{equation}

A hypersurface is called parallel \cite{Asperti} if $\tilde{\nabla}h=0$.

\begin{th}
Let $G$ be an hypersurface of a generalized f.pk-space forms $%
M^{2n+s}(F_{1},F_{2},{\cal F})$, tangent to the structure vector fields with
$F_{2}\not=0$, then $G$ is not parallel.
\end{th}

Proof. Let $G$ be a parallel hypersurface of a generalized f.pk-space forms $%
M^{2n+s}(F_{1},F_{2},{\cal F})$, tangent to the structure vector fields. Let
$h$ be the second fundamental form of $G$. Consider $N$ be the unit normal
of $G$ in $M$ and $W=-\varphi N$. Since $g(\xi _{i},N)=\eta _{i}(N)=0$ for
all $i$, therefore
\[
g(W,W)=g(\varphi N,\varphi N)=1,
\]%
\[
g(N,W)=g(N,\varphi N)=0,
\]%
\[
\varphi W=N.
\]%
Let $X\in \Gamma (TG)$, then we have
\begin{equation}
\varphi X=TX+w(X)N  \label{eq-phi}
\end{equation}%
where $w$ and $T$ are tensor fields on $G$ of type $(0,1)$ and $(1,1)$,
respectively and $TX$ is also the tangent part of $\varphi X$.

Set $w\not=0$, by (\ref{eq-phi}), we have $w(X)=g(X,W)$.

Now, using (\ref{eq-fpk}), (\ref{eq-phi}) and Codazzi equation, we have
\begin{eqnarray*}
0 &=&\tilde{\nabla}_{X}h(Y,Z)-\tilde{\nabla}_{Y}h(X,Z)=\left( \tilde{R}%
(X,Y)Z\right) ^{\bot } \\
&=&F_{2}\left( g(Z,\varphi Y)\varphi X-g(Z,\varphi X)\varphi Y+2g(X,\varphi
Y)\varphi Z\right) ^{\bot }
\end{eqnarray*}%
Since $F_{2}\not=0$, we obtain
\[
\left( g(Z,\varphi Y)w(X)-g(Z,\varphi X)w(Y)+2g(X,\varphi Y)w(Z)\right) N=0.
\]%
Putting $Z=W$ in above equation, we get
\[
2g(X,TY)=0,
\]%
which implies that $TY=0$, so $\dim \varphi (T_{p}\left( G\right) )=1$ for
all $p\in G$.

Since $T_{p}\left( M\right) =T_{p}\left( G\right) \oplus \left( T_{p}\left(
G\right) \right) ^{\bot }$ and $rank\varphi =2n$, therefore we obtain
\[
2n-1\leq \dim \varphi \left( T_{p}\left( G\right) \right) ^{\bot }\leq 2n,
\]%
which is impossible since $n>1$.

\begin{rem-new}
This result is true for $S$-space form \cite{Belkhelfa} and generalized
Sasakian-space form \cite{Gherib}.
\end{rem-new}

\begin{th}
There are semi-parallel hypersurfaces, tangent to the structure vector
fields in a generalized f.pk-space forms $M^{2n+s}(F_{1},F_{2},{\cal F})$
with assumption that ${\cal F}$ has atleast one non-zero function.
\end{th}

\noindent {\bf Proof.} Let $G$ be not a semi-parallel hypersurfaces, tangent
to the structure vector fields in a generalized f.pk-space forms $%
M^{2n+s}(F_{1},F_{2},{\cal F})$. It means that
\[
\tilde{R}.h(X,Y,Z,W)=-h(\tilde{R}(X,Y)Z,W)-h(Z,\tilde{R}(X,Y)W)\not=0,
\]%
where $h$ is the second fundamental form of $G$, which is symmetric $(0,2)$%
-tensor field on $G$. By using the same argument of Theorem \ref{th-sec}, we
get
\[
H(X,Z)\not=g(X,Z)\sum\limits_{\gamma =1}^{s}H(\xi _{j},\xi _{\gamma
})-\sum\limits_{\alpha ,\beta ,\gamma =1}^{s}\eta ^{\alpha }(X)\eta ^{\beta
}(Z)H(\xi _{j},\xi _{\gamma })+\sum\limits_{\alpha ,\beta ,\gamma
=1}^{s}\eta ^{\alpha }(X)\eta ^{j}(Z)H(\xi _{\beta },\xi _{\gamma }),
\]%
which is impossible.

\begin{rem-new}
Even for $S$-space form \cite{Belkhelfa}, this result is not true.
\end{rem-new}

\bigskip

\noindent Punam Gupta \\
Department of Mathematics and Statistics
\newline Dr. Harisingh Gour University\newline Sagar-470 003, Madhya
Pradesh\newline India\newline Email:{\em punam\_2101\makeatletter
@\makeatother yahoo.co.in}\\

\medskip

\noindent Sanjay Kumar Singh\\
Department of Mathematics \newline Indian Institute of Science
Education \& Research\newline Bhopal-462 066, Madhya Pradesh\newline
India\newline Email:{\em sanjayks\makeatletter @\makeatother
iiserb.ac.in}

\begin{thebibliography}{99}
\bibitem{Arslan} K. Arslan, R. Ezenta\c{s}, C. Murathan and C. \"{O}zg\"{u}%
r, On pseudo Ricci-symmetric manifolds, Balkan J. Geom. Appl. 6 (2001), no.
2, 1-5.

\bibitem{Asperti} C. Asperti, A. Lobos, F. Mercuri, Pseudo-parallel
immersions of a space forms, Adv. Geom. 2(2002), no. 1, 57-71.

\bibitem{Bejan} C.L. Bejan and M. Crasmareanu, Parallel second-order tensors
on Vaisman manifolds, Int. J. Geom. Methods Mod. Phys. 14 (2017), no. 2,
1750023, 8 pp.

\bibitem{Belkhelfa} M. Belkhelfa and F. Mahi, Second order parallel tensors
on ${\cal S}$-manifolds and semi-parallel hypersurfaces of ${\cal S}$-space
forms, Ukra\"{\i}n. Mat. Zh. 71 (2019), no. 10, 1422-1429.

\bibitem{Blair} D.E. Blair, Geometry of manifolds with structural group $%
U(n)\times O(s)$, J. Differ. Geom. 4(1970), 155-167.

\bibitem{Chaki} M.C. Chaki, On pseudo-Ricci symmetric manifolds. Bulg. J.
Phys. 15 (1988), 525-531.

\bibitem{Das} L. Das, Second order parallel tensors on $\alpha $-Sasakian
manifold, Acta Math. Acad. Paedagog. Nyh\'{a}zi. (N.S.) 23(2007), no. 1,
65-69.

\bibitem{Deszcz} R. Deszcz, On Ricci-pseudosymmetric warped products,
Demonstratio Math., 22 (1989), 1053-1065.

\bibitem{Duggal} K. Duggal, S. Ianus and A.M. Pastore, Maps interchanging $f$%
-structures and their harmonicity, Acta Appl. Math. 67(2001), no. 1, 91-115.

\bibitem{Eisenhart-23} L.P. Eisenhart, Symmetric tensors of the second order
whose first covariant derivatives are zero, Trans. Amer. Math. Soc. 25
(1923), 297-306.

\bibitem{Fal-Pas} M. Falcitelli and A. M. Pastore, Generalized globally
framed $f$-space-forms, Bull. Math. Soc. Sci. Math. Roumanie, 52(2009), no.
3, 291-305.

\bibitem{Gherib} F. Gherib, M. Belkhelfa, Parallel submanifolds of
generalized Sasakian space forms, Bull. Transilv. Univ. Bra\c{s}ov Ser. III,
51(2009), no. 2, 185-191.

\bibitem{Gherib-1} F. Gherib, and M. Belkhelfa, Second order parallel
tensors on generalized Sasakian space forms and semi parallel hypersurfaces
in Sasakian space forms, Beitr\"{a}ge Algebra Geom. 51 (2010), no. 1, 1-7.

\bibitem{Gold} S.I. Goldberg and K. Yano, Globally framed $f$-manifolds,
Ill. J. Math. 15(1971), 456-474.

\bibitem{Gupta-1} P. Gupta, Index of pseudo-projectively symmetric
semi-Riemannian manifolds, Carpathian Math. Publ. 7(2015), no. 1, 57-65.

\bibitem{Kobayashi} M. Kobayashi and S. Tsuhiya, Invariant submanifolds of
an $f$-manifolds with complemented frames, Kodai Math. Semi. Rep, 24(1972),
430-450.

\bibitem{Levy-26} H. Levy, Symmetric tensors of the second order whose
covariant derivatives vanish, Annals of Maths. 27(1926), 91-98.

\bibitem{Mondal} A.K. Mondal, U.C. De and C. \"{O}zg\"{u}r, Second order
parallel tensors on $\left( \kappa ,\mu \right) $-contact metric manifolds,
An. \c{S}tiin\c{t}. Univ. "Ovidius\textquotedblright\ Constan\c{t}a Ser.
Mat. 18 (2010), no. 1, 229-238.

\bibitem{Sharma} R. Sharma, Second order parallel tensor in real and complex
space forms, Internat. J. Math.Math. Sci. 12(1989), no. 4, 787-790.

\bibitem{Ramesh} R. Sharma, Second order parallel tensors on contact
manifolds, Algebras Groups Geom. 7(1990), no. 2, 145-152.

\bibitem{Sharma-1} R. Sharma, Second order parallel tensors on contact
manifolds. II, C. R. Math. Rep. Acad. Sci. Canada 13 (1991), no. 6, 259-264.

\bibitem{Gupta} M.M. Tripathi, P. Gupta and J.-S. Kim, Index of
quasiconformally symmetric semi-Riemannian manifolds, Int. J. Math. Math.
Sci. 2012(2012), 14 pp.

\bibitem{Yano} K. Yano, On a structure defined by a tensor field $f$
satisfying $f^{3}+f=0$, Tensor N. S. 14(1963), 99-109.
\end{thebibliography}
\end{document}